\documentclass[oneside,english]{amsart}
\usepackage[T1]{fontenc}
\usepackage[latin9]{inputenc}
\usepackage{amsthm}
\usepackage{amssymb}

\makeatletter
\numberwithin{equation}{section}
\numberwithin{figure}{section}
\theoremstyle{plain}
\newtheorem{thm}{\protect\theoremname}

\makeatother

\usepackage{babel}
\providecommand{\theoremname}{Theorem}

\begin{document}

\title{Zeros of Ramanujan Type Entire Functions}

\author{Ruiming Zhang}
\begin{abstract}
In this work we establish some polynomials and entire functions have
only real zeros. These polynomials generalize q-Laguerre polynomials
$L_{n}^{(\alpha)}(x;q)$, while the entire functions are generalizations
of Ramanujan's entire function $A_{q}(z)$, q-Bessel functions $J_{\nu}^{(2)}(z;q)$,
$J_{\nu}^{(3)}(z;q)$ and confluent basic hypergeometric series.
\end{abstract}


\curraddr{College of Science\\
Northwest A\&F University\\
Yangling, Shaanxi 712100\\
P. R. China.}


\thanks{This research is partially supported by National Natural Science
Foundation of China, grant No. 11371294 and Northwest A\&F University.}

\maketitle

\section{Introduction}

In \cite{Ismail2} Ismail and Zhang defined an entire function 
\begin{equation}
A_{q}^{\left(\alpha\right)}\left(a;z\right)=\sum_{n=0}^{\infty}\frac{\left(a;q\right)_{n}q^{\alpha n^{2}}z^{n}}{\left(q;q\right)_{n}},\label{eq:1.1}
\end{equation}
where $\alpha>0$, $0<q<1$, and for any positive integer $m\in\mathbb{N}$
and all complex numbers $a,a_{1},\dots,a_{m}\,n\in\mathbb{C}$, \cite{Andrews1,Ismail1}
\begin{equation}
(a;q)_{\infty}=\prod_{k=0}^{\infty}\left(1-aq^{k}\right),\ (a;q)_{n}=\frac{\left(a;q\right)_{\infty}}{\left(aq^{n};q\right)_{\infty}},\ (a_{1},\dots,a_{m};q)_{n}=\prod_{j=1}^{m}(a_{j};q)_{n}.\label{eq:1.2}
\end{equation}
Clearly, $A_{q}^{\left(\alpha\right)}\left(a;z\right)$ is a kind
of entire function treated in Lemma 14.1.4 of \cite{Ismail1}, and
by the lemma it must have infinitely many zeros if $aq^{n}\neq1$
for any nonnegative integers $n\in\mathbb{N}_{0}$. Furthermore, in
\cite{Hayman} Hayman proved an asymptotic expansion for the $n$-th
zero of this class of entire function under the general condition
that the parameter $q$ is strictly inside the unit disk $|q|<1$. 

We observe that

\begin{equation}
A_{q}^{\left(1/2\right)}\left(q^{-n};z\right)=\sum_{k=0}^{\infty}\frac{\left(q^{-n};q\right)_{k}q^{k^{2}/2}z^{k}}{\left(q;q\right)_{k}}=\left(q;q\right)_{n}S_{n}\left(zq^{1/2-n};q\right),\label{eq:1.3}
\end{equation}
 and

\begin{equation}
A_{q}^{\left(1\right)}\left(0;z\right)=\sum_{n=0}^{\infty}\frac{q^{n^{2}}z^{n}}{(q;q)_{n}}=A_{q}\left(-z\right),\ A_{q}^{\left(1\right)}\left(q;z\right)=\sum_{n=0}^{\infty}q^{n^{2}}z^{n},\label{eq:1.4}
\end{equation}
where $A_{q}(z)$ and $S_{n}(z;q)$ are the Ramanujan entire function
and Stieltjes-Wigert polynomial respectively, see \cite{Ismail1},
for example. Since $A_{q}^{\left(\alpha\right)}\left(a;z\right)$
generalizes both $A_{q}(z)$ and $S_{n}(z;q)$, and it is well-known
that both of them have all real positive zeros. Thus it is natural
to ask under what conditions the zeros of the entire function $A_{q}^{\left(\alpha\right)}\left(a;z\right)$
are all real. In this work we shall present a partially answer to
this question and more.

\section{Preliminaries}

In the proofs we need the Vitali's theorem \cite{Titchmarsh1} and
Hurwitz's theorem \cite{Ahlfors1}. For our convenience we list them
here. The first is the Vitali's theorem:
\begin{thm}
\label{thm:Vitali} Let $\left\{ f_{n}(z)\right\} $ be a sequence
of functions analytic in a domain $D$ and assume that $f_{n}(z)\to f(z)$
point-wise in $D$. Then $f_{n}(z)\to f(z)$ uniformly in any subdomain
bounded by a contour $C$, provided that $C$ is contained in $D$.
\end{thm}
Here is the Hurwitz's theorm:
\begin{thm}
\label{thm:Hurwitz} If the functions $f_{n}(z)$ are analytic and
$\neq0$ in a region $\Omega$, and if $f_{n}(z)$ converges to $f(z)$,
uniformly on every compact subset of $\Omega$, then $f(z)$ is either
identically zero or never equal to zero in $\Omega$.
\end{thm}
It is known that if an entire function has a finite order that it
is not a positive integer, then $f(z)$ has infinitely many zeros,
\cite{Boas,Ismail1}. Given an entire function of order $0$ with
$f(0)\neq0$, if $-z_{k},\ k\in\mathbb{N}$ are all the roots, then
by Hadamard's canonical representation of entire functions we have
\begin{equation}
f(z)=f(0)\prod_{k=1}^{\infty}\left(1+\frac{z}{z_{k}}\right),\label{eq:2.2}
\end{equation}
where ${\displaystyle \sum_{k=1}^{\infty}}|z_{k}|^{-1}<\infty.$ A
real entire function $f(z)$ is of Laguerre-P\'{o}lya class if \cite{Dimitrov1,Karlin}
\begin{equation}
f(z)=cz^{m}e^{-\alpha z^{2}+\beta z}\prod_{k=1}^{\infty}\left(1+\frac{z}{z_{k}}\right)e^{-z/z_{k}},\label{eq:2.3}
\end{equation}
where $c,\beta,\,z_{k}\in\mathbb{R}$, $\alpha\ge0$, $m\in\mathbb{Z}^{+}$,
and ${\displaystyle \sum_{k=1}^{\infty}}z_{k}^{-2}<\infty$. Clearly,
given a real entire function $f(z)$ of order $0$ that has nonnegative
Taylor coefficients and satisfies $f(0)\neq0$, if it is also in the
Laguerre-P\'{o}lya class then it must have the factorization (\ref{eq:2.2})
with $z_{k}>0$. 

A real sequence $\left\{ a_{n}\right\} _{n=0}^{\infty}$ is called
a P\'{o}lya frequence (PF) sequence if the infinite matrix $\left(a_{j-i}\right)_{i,j=0}^{\infty}$
is totally positive, i.e. all its minors are nonnegative, where we
follow the usual convention that $a_{n}=0$ if $n<0$. \cite{Aissen1,Driver1,Karlin}
\begin{thm}
\label{thm:aesw} The sequence $\left\{ b_{k}\right\} _{k=0}^{\infty}$
is a PF sequence if and only if 
\begin{equation}
cz^{m}e^{\gamma z}\prod_{k=1}^{\infty}\frac{1+\alpha_{k}z}{1-\beta_{k}z}=\sum_{k=0}^{\infty}b_{k}z^{k},\label{eq:2.5}
\end{equation}
where $c\ge0,\ \gamma\ge0,\ \alpha_{k}\ge0,\ \beta_{k}\ge0$, $m\in\mathbb{Z}^{+}$,
and ${\displaystyle \sum_{k=1}^{\infty}}\left(\alpha_{k}+\beta_{k}\right)<\infty$. 
\end{thm}
In the case $b_{k}=0$ for $k=n+1,\dots$, then the right hand series
becomes a polynomial of degree at most $n$, then we have the following,
\cite{Aissen1,Driver1,Karlin,Pitman}
\begin{thm}
\label{thm:pf}The PF sequences have the following properties:
\begin{enumerate}
\item Let $b_{k}\ge0,\ k\ge0$. Then, the sequence $\left\{ b_{k}\right\} _{k=0}^{n}$
is a PF sequence if and only if the polynomial ${\displaystyle \sum_{k=0}^{n}}b_{k}x^{k}$
has all nonpositive zeros.
\item Let $a_{k},\,b_{k}\ge0,\ k\ge0$. If $\left\{ a_{k}\right\} _{k=0}^{m}$
and $\left\{ b_{k}\right\} _{k=0}^{n}$ are PF sequences, then so
is $\left\{ a_{k}b_{k}\right\} _{k\ge0}$.
\item Let $b_{0},\,b_{1},\,\dots,\,b_{n}\ge0$ and $\left\{ b_{k}\right\} _{k=0}^{n}$
be a PF sequence, then so is $\left\{ b_{k}/k!\right\} _{k=0}^{n}$
.
\item Let $\left\{ a_{k}\right\} _{k=0}^{m}$ and $\left\{ b_{k}\right\} _{k=0}^{n}$
be PF sequences such that $a_{k},\,b_{k}\ge0,\ k\ge0$, then $\left\{ k!a_{k}b_{k}\right\} _{k\ge0}$
is also a PF sequence.
\end{enumerate}
\end{thm}
It was proved in that \cite{Carnicer1,Dimitrov1} 
\begin{equation}
\sum_{n=0}^{\infty}\frac{q^{n^{2}}}{n!}x^{n}\label{eq:2.6}
\end{equation}
is in Laguerre-P\'{o}lya class for $q\in(-1,1)$. For $0<q<1$. This
function is clearly a real entire function of order $0$ with positive
Taylor coefficients, hence the sequence $\left\{ q^{n^{2}}/n!\right\} _{n=0}^{\infty}$
is a PF sequence by Theorem \ref{thm:aesw}.

For $a\ge0$ and $0<q<1$, from the $q$-Binomial theorem we have
\cite{Andrews1,Ismail1}
\begin{equation}
\frac{(-az;q)_{\infty}}{(bz;q)_{\infty}}=\sum_{n=0}^{\infty}\frac{(-a;q)_{n}}{(q;q)_{n}}z^{n},\label{eq:2.7}
\end{equation}
and
\begin{equation}
\left(-q^{n}z;q\right)_{n}=\sum_{k=0}^{n}\frac{\left(q^{-n};q\right)_{k}\left(-z\right)^{k}}{(q;q)_{k}}.\label{eq:2.8}
\end{equation}
Then by Theorem \ref{thm:aesw} both$\left\{ (-a;q)_{k}/(q;q)_{k}\right\} _{k=0}^{\infty}$
and $\left\{ \left(-1\right)^{k}\left(q^{-n};q\right)_{k}/(q;q)_{k}\right\} _{k=0}^{\infty}$
are PF sequences.

The $q$-Bessel functions $J_{\nu}^{(1)}(z;q)$, $J^{(2)}(z;q)$ and
$J^{(3)}(z;q)$ are defined by \cite{Ismail1}

\begin{equation}
\frac{2^{\nu}(q;q)_{\infty}J_{\nu}^{(1)}(z;q)}{\left(q^{\nu+1};q\right)_{\infty}z^{\nu}}=\sum_{n=0}^{\infty}\frac{\left(-z^{2}/4\right)^{n}}{\left(q,q^{\nu+1};q\right)_{n}},\label{eq:2.9}
\end{equation}
\begin{equation}
\frac{2^{\nu}(q;q)_{\infty}J_{\nu}^{(2)}(z;q)}{\left(q^{\nu+1};q\right)_{\infty}z^{\nu}}=\sum_{n=0}^{\infty}\frac{q^{n^{2}}\left(-z^{2}q^{\nu}/4\right)^{n}}{\left(q,q^{\nu+1};q\right)_{n}},\label{eq:2.10}
\end{equation}
 and
\begin{eqnarray}
\frac{2^{\nu}(q;q)_{\infty}J_{\nu}^{(3)}(z;q)}{\left(q^{\nu+1};q\right)_{\infty}z^{\nu}} & = & \sum_{n=0}^{\infty}\frac{q^{\binom{n+1}{2}}\left(-z^{2}/4\right)^{n}}{\left(q,q^{\nu+1};q\right)_{n}}\label{eq:2.11}
\end{eqnarray}
respectively.

For $\nu>-1$ and $0<q<1$, it is known that $J^{(2)}(z;q)$ is an
entire function of order $0$ such that 
\begin{equation}
\frac{2^{\nu}(q;q)_{\infty}J_{\nu}^{(2)}(z;q)}{\left(q^{\nu+1};q\right)_{\infty}z^{\nu}}=\prod_{n=0}^{\infty}\left(1-\frac{z^{2}}{j_{\nu,n}(q)^{2}}\right),\label{eq:2.13}
\end{equation}
 where 
\begin{equation}
0<j_{\nu,1}(q)<j_{\nu,2}(q)<\dots\label{eq:2.14}
\end{equation}
 are the positive zeros of $J^{(2)}(z;q)$ satisfying 
\begin{equation}
\sum_{n=1}^{\infty}\frac{1}{j_{\nu,n}(q)^{2}}<\infty.\label{eq:2.15}
\end{equation}
For $\nu>-1$ and $0<q<1$, from (\ref{eq:2.13}) and the relation
\cite{Ismail1} 
\begin{equation}
J_{\nu}^{(1)}(z;q)=\frac{J_{\nu}^{(2)}(z;q)}{\left(-z^{2}/4;q\right)_{\infty}}\label{eq:2.16}
\end{equation}
we get 
\begin{equation}
\sum_{n=0}^{\infty}\frac{z^{n}}{\left(q,q^{\nu+1};q\right)_{n}}=\frac{{\displaystyle \prod_{n=0}^{\infty}}\left(1+z/j_{\nu,n}(q)^{2}\right)}{\left(z/4;q\right)_{\infty}}.\label{eq:2.17}
\end{equation}
From (\ref{eq:2.10}), (\ref{eq:2.13}) and (\ref{eq:2.17}) we see
that both $\left\{ 1/\left(q,q^{\nu+1};q\right)_{k}\right\} _{k=0}^{\infty}$
and $\left\{ q^{k^{2}}/\left(q,q^{\nu+1};q\right)_{k}\right\} _{k=0}^{\infty}$
are PF sequences by Theorem \ref{thm:aesw}. 

It is worth noticing that these $q$-Bessel functions are much more
than merely being $q$-analogues to their classical counterpart 
\begin{equation}
J_{\nu}(z)=\frac{\left(z/2\right)^{\nu}}{\Gamma(\nu+1)}\sum_{n=0}^{\infty}\frac{\left(-z^{2}/4\right)^{n}}{n!(\nu+1)_{n}},\label{eq:2.18}
\end{equation}
where
\begin{equation}
(a)_{n}=\frac{\Gamma(a+n)}{\Gamma(a)},\quad a,n\in\mathbb{Z}\label{eq:2.19}
\end{equation}
with $-a-n\notin\mathbb{N}_{0}$ and $\Gamma(z)$ is the Euler's Gamma
function. They have many interesting properties that the classical
$J_{\nu}(z)$ does not possess. Here we only mention two of them.
First of them is that $J_{\nu}^{(2)}\left(2iq^{-n/2};q\right),\ n\in\mathbb{N}_{0}$
can be evaluated explicitly, see \cite{Ismail3}. The second is that
there are many symmetries among themselves. 

Let
\begin{equation}
j_{\nu}^{(k)}(z;q)=\frac{2^{\nu}(q;q)_{\infty}J_{\nu}^{(k)}(z;q)}{\left(q^{\nu+1};q\right)_{\infty}z^{\nu}},\quad k=1,2,3.\label{eq:2.20}
\end{equation}
Then for $q>1$ and $|z|<1$, we have
\begin{equation}
j_{\nu}^{(1)}\left(z;\frac{1}{q}\right)=j_{\nu}^{(2)}\left(\sqrt{q}z;q\right),\ j_{\nu}^{(2)}\left(z;\frac{1}{q}\right)=j_{\nu}^{(1)}\left(\sqrt{q}z;q\right),\ j_{2}^{(3)}\left(z;\frac{1}{q}\right)=j_{\nu}^{(3)}\left(q^{\frac{\nu}{2}}z;q\right).\label{eq:2.21}
\end{equation}
The above relations are very typical among basic hypergeometric series,
but it is seldom explored in the literature.

\section{Main Results}

We follow the usual convention that an empty sum is zero but an empty
product is 1.

\begin{thm}
\label{thm:poly} Let $n\in\mathbb{N}$, $0<q<1$, and $\alpha\ge0$.
\begin{enumerate}
\item The polynomial 
\begin{equation}
A_{q}^{(\alpha)}\left(q^{-n};x\right)=\sum_{k=0}^{n}\frac{\left(q^{-n};q\right)_{k}q^{\alpha k^{2}}x^{k}}{\left(q;q\right)_{k}}\label{eq:3.2}
\end{equation}
 has all positive zeros. 
\item For $m,\ell\in\mathbb{N}_{0}$, $n_{j}\in\mathbb{N},\ 0<q_{j}<1,\ 1\le j\le m$,
and $\beta_{r}>0,\ 0<q_{r}<1,\ 1\le r\le\ell$, the polynomial
\begin{equation}
\sum_{k=0}^{\min\left\{ n_{j}\vert1\le j\le m\right\} }\prod_{j=1}^{m}\frac{\left(q_{j}^{-n_{j}};q_{j}\right)_{k}}{(q_{j};q_{j})_{k}}\frac{q^{\alpha k^{2}}\left((-1)^{m}x\right)^{k}}{{\displaystyle \prod_{r=1}^{\ell}}\left(q_{r},q_{r}^{\beta_{r}};q_{r}\right)_{k}}\label{eq:3.3}
\end{equation}
has all negative zeros. 
\item For $m,\ell\in\mathbb{N}_{0}$, $n_{j}\in\mathbb{N},\ 1\le j\le m$,
and $\beta_{r}>0,\ 1\le r\le\ell$, the polynomial
\begin{equation}
\sum_{k=0}^{\min\left\{ n_{j}\vert1\le j\le m\right\} }\frac{\prod_{j=1}^{m}(-n_{j})_{k}}{\prod_{r=1}^{\ell}\left(\beta_{r}\right)_{k}}\frac{\left((-1)^{m}x\right)^{k}}{(k!)^{m+\ell}}\label{eq:3.4}
\end{equation}
has all negative zeros.
\end{enumerate}
\end{thm}
The following results are for the entire functions in the first parameter
range:
\begin{thm}
\label{thm:func1} Let $0<q<1$ and $\alpha>0$. Then,
\begin{enumerate}
\item For $a\ge0$, the entire function 
\begin{equation}
A_{q}^{\left(\alpha\right)}\left(-a;z\right)=\sum_{n=0}^{\infty}\frac{\left(-a;q\right)_{n}q^{\alpha n^{2}}z^{n}}{\left(q;q\right)_{n}}\label{eq:3.5}
\end{equation}
has infinitely many zeros and all of them are negative. 
\item For $m,\ell\in\mathbb{N}_{0}$, $a_{j}\ge0,\ 0<q_{j}<1,\ 1\le j\le m$,
and $\beta_{r}>0,\ 0<q_{r}<1,\ 1\le r\le\ell$, the entire function
\begin{equation}
\sum_{k=0}^{\infty}\prod_{j=1}^{m}\frac{\left(-a_{j};q_{j}\right)_{k}}{(q_{j};q_{j})_{k}}\frac{q^{\alpha k^{2}}z^{k}}{{\displaystyle \prod_{r=1}^{\ell}}\left(q_{r},q_{r}^{\beta_{r}};q_{r}\right)_{k}}\label{eq:3.6}
\end{equation}
has infinitely many zeros and all of them are negative. 
\item For $m\in\mathbb{N}_{0},\ \ell\in\mathbb{N}$ and $\beta_{r}\ge1,\ 1\le r\le\ell$,
the entire function
\begin{equation}
\sum_{k=0}^{\infty}\frac{z^{k}}{\left(k!\right)^{m+\ell}{\displaystyle \prod_{r=1}^{\ell}\left(\beta_{r}\right)_{k}}}\label{eq:3.7}
\end{equation}
has infinitely many zeros and all of them are negative. 
\end{enumerate}
\end{thm}
For $r,\,s\in\mathbb{N}_{0}$ and $0<q<1$, let $a_{1},\dots,\,a_{r},\ b_{1},\dots,\,b_{s}\in\mathbb{C}$
\begin{equation}
_{r}A_{s}^{(\alpha)}(a_{1},\dots a_{r};\ b_{1},\dots,\ b_{s}\ ;\,q;\,z)=\sum_{n=0}^{\infty}\frac{(a_{1},\dots,\,a_{r};\,q)_{n}}{\left(b_{1},\dots,b_{s};\,q\right)_{n}}q^{\alpha n^{2}}z^{n},\label{eq:3.8}
\end{equation}
 then it is clear that for $\alpha>0$ and $0<q<1$ we have 
\begin{equation}
A_{q}^{(\alpha)}(a;z)=_{1}A_{1}^{(\alpha)}(a,\,q;z),\label{eq:3.9}
\end{equation}
 and for $z\in\mathbb{C}$ and $s+1>r$ we have 
\begin{equation}
\begin{aligned} & _{r}A_{s}^{\left((s+1-r)/2\right)}\left(a_{1},\dots a_{r};\ q,\,b_{1},\dots,\ b_{s}\ ;\,q;\,(-1/\sqrt{q})^{s+1-r}z\right)\\
 & ={}_{r}\phi_{s}\left(\begin{array}{cc}
\begin{array}{c}
a_{1},\dots,\,a_{r}\\
b_{1},\dots,\,b_{s}
\end{array} & \bigg|q,z\end{array}\right),
\end{aligned}
\end{equation}
where the basic hypergeometric series $_{r}\phi_{s}$ is defined by
\cite{Andrews1,Ismail1}
\begin{equation}
\begin{aligned}_{r}\phi_{s}\left(\begin{array}{cc}
\begin{array}{c}
a_{1},\dots,\,a_{r}\\
b_{1},\dots,\,b_{s}
\end{array} & \bigg|q,z\end{array}\right) & =\sum_{n=0}^{\infty}\frac{(a_{1},\dots,\,a_{r};q)_{n}}{(q,b_{1},\dots,\,b_{s};q)_{n}}\left(-q^{(n-1)/2}\right)^{n(s+1-r)}z^{n}.\end{aligned}
\label{eq:61}
\end{equation}
 Finally, here is the result for entire functions in another parameter
range:
\begin{thm}
\label{thm:func2} For $r,\,s\in\mathbb{N}_{0}$, let $\alpha>0,\ 1>a_{j},\,b_{k}>0,\ 1\le j\le r,\ 1\le k\le s$
such that 
\begin{equation}
\prod_{j=1}^{r}\left(1-a_{j}\right)\prod_{k=1}^{s}\left(1-b_{k}q\right)\ge4q^{2\alpha},\label{eq:3.10}
\end{equation}
then the entire function $_{r}A_{s}^{(\alpha)}(a_{1},\dots a_{r};\ b_{1},\dots,\ b_{s}\ ;\,q;\,z)$
has all negative zeros. 
\end{thm}

\section{Proofs}

\subsection{Proof of Theorem \ref{thm:poly}}

For $n\in\mathbb{N}$, $0<q<1$ and $\alpha>0$, since the sequences
$\left\{ \left(q^{-n};q\right)_{k}\left(-1\right)^{k}/(q;q)_{k}\right\} _{k=0}^{n}$
and $\left\{ q^{\alpha k^{2}}/k!\right\} _{k=0}^{n}$are PF sequences,
then apply property 4 of Theorem \ref{thm:pf} to get the PF sequence
$\left\{ \left(q^{-n};q\right)_{k}\left(-1\right)^{k}q^{\alpha k^{2}}/(q;q)_{k}\right\} _{k=0}^{n}$.
Then by property 1 of Theorem \ref{thm:pf} we know the polynomial
\[
\sum_{k=0}^{n}\frac{\left(q^{-n};q\right)_{k}\left(-1\right)^{k}q^{\alpha k^{2}}z^{k}}{(q;q)_{k}}=A_{q}^{(\alpha)}\left(q^{-n};-z\right)
\]
has all negative zeros, which is equivalent to that $A_{q}^{(\alpha)}\left(q^{-n};z\right)$
has all positive zeros.

For $m,\ell\in\mathbb{N}$, $n_{j}\in\mathbb{N},\ 0<q_{j}<1,\ 1\le j\le m$,
and $\nu_{r}>-1,\ 0<q_{r}<1,\ 1\le r\le\ell$, first by property 2,
then by property 4 of Theorem \ref{thm:pf} we know that 
\[
\left\{ \prod_{j=1}^{m}\frac{\left(q_{j}^{-n_{j}};q_{j}\right)_{k}}{(q_{j};q_{j})_{k}}\frac{\left(-1\right)^{km}}{{\displaystyle \prod_{r=1}^{\ell}}\left(q_{r},q_{r}^{\nu_{r}+1};q_{r}\right)_{k}}\right\} _{k\ge0}
\]
 and
\[
\left\{ \prod_{j=1}^{m}\frac{\left(q_{j}^{-n_{j}};q_{j}\right)_{k}}{(q_{j};q_{j})_{k}}\frac{\left(-1\right)^{km}q^{\alpha k^{2}}}{{\displaystyle \prod_{r=1}^{\ell}}\left(q_{r},q_{r}^{\nu_{r}+1};q_{r}\right)_{k}}\right\} _{k\ge0}
\]
are all PF sequences. Then for $\alpha\ge0$ the polynomial (\ref{eq:3.3})
has all negative zeros. 

For $m,\ell\in\mathbb{N}$, take 
\[
q_{j}=q,\ 1\le j\le m,\quad q_{r}=q,\ 1\le r\le\ell,\quad z=\left(1-q\right)^{2\ell}x,
\]
 and let $q\uparrow1$ in 
\[
\sum_{k=0}^{\min\left\{ n_{j}\vert1\le j\le m\right\} }\prod_{j=1}^{m}\prod_{s=1}^{k}\frac{\left(1-q^{s-1-n_{j}}\right)}{\left(1-q^{s}\right)}\frac{\left(-1\right)^{km}q^{\alpha k^{2}}\left(1-q\right)^{2k\ell}x^{k}}{{\displaystyle \prod_{r=1}^{\ell}}{\displaystyle \prod_{t=1}^{k}}\left(1-q^{t}\right)\left(1-q^{\nu_{r}+t}\right)}
\]
 to get
\[
\begin{aligned} & \sum_{k=0}^{\min\left\{ n_{j}\vert1\le j\le m\right\} }\prod_{j=1}^{m}\prod_{s=1}^{k}\frac{s-1-n_{j}}{s}\frac{\left(-1\right)^{km}x^{k}}{{\displaystyle \prod_{r=1}^{\ell}}{\displaystyle \prod_{t=1}^{k}}t(\nu_{r}+t)}\\
 & =\sum_{k=0}^{\min\left\{ n_{j}\vert1\le j\le m\right\} }\frac{\prod_{j=1}^{m}(-n_{j})_{k}}{{\displaystyle \prod_{r=1}^{\ell}}\left(\nu_{r}+1\right)_{k}}\frac{\left(-1\right)^{km}x^{k}}{(k!)^{m+\ell}}
\end{aligned}
\]
for each $z\in\mathbb{C}$. Then this limit is also uniform on any
compact subset of $\mathbb{C}$ by Vitali's theorem \ref{thm:Vitali}.
Thus the polynomial (\ref{eq:3.4}) has all negative zeros by applying
Hurwitz's theorem \ref{thm:Hurwitz}.

\subsection{Proof of Theorem \ref{thm:func1}}

For $n\in\mathbb{N}$, $0<q<1$, $a\ge0$ and $\alpha>0$, since the
sequences $\left\{ \left(-a;q\right)_{k}/(q;q)_{k}\right\} _{k=0}^{n}$
and $\left\{ q^{\alpha k^{2}}/k!\right\} _{k=0}^{n}$are PF sequences,
then apply property 4 of Theorem \ref{thm:pf} to get the PF sequence
$\left\{ \left(-a;q\right)_{k}q^{\alpha k^{2}}/(q;q)_{k}\right\} _{k=0}^{n}$.
Then by property 1 of Theorem \ref{thm:pf} we know the polynomial
\[
\sum_{k=0}^{n}\frac{\left(-a;q\right)_{k}q^{\alpha k^{2}}}{(q;q)_{k}}z^{k}
\]
has all negative zeros. For $0<q<1$, $a\ge0$ and $\alpha>0$ and
for each $n\in\mathbb{N}$ and $z\in\mathbb{C}$ we have, 
\[
\left|\sum_{k=0}^{n}\frac{\left(-a;q\right)_{k}q^{\alpha k^{2}}}{(q;q)_{k}}z^{k}\right|\le\sum_{k=0}^{\infty}\frac{\left(-a;q\right)_{k}q^{\alpha k^{2}}}{(q;q)_{k}}\left|z\right|^{k}<\infty.
\]
Now we first apply Vitali's theorem, then apply the Hurwitz's theorem
to show the entire function $A_{q}^{(\alpha)}(-a;z)$ has no zeros
outside the set $(-\infty,0)$. Since by Lemma 14.1.4 of \cite{Ismail1}
we know that $A_{q}^{(\alpha)}(-a;z)$ has infinitely many zeros,
then we have proved that $A_{q}^{(\alpha)}\left(-a;z\right)$ has
infinitely many zeros and all of them are negative.

For $m,\ell\in\mathbb{N}_{0},\ \alpha>0,\ 0<q<1$, $a_{j}\ge0,\ 0<q_{j}<1,\ 1\le j\le m$,
and $\nu_{r}>-1,\ 0<q_{r}<1,\ 1\le r\le\ell$, similar to the polynomial
case, we first apply property 2, then apply property 4 of Theorem
\ref{thm:pf} we know that for each $n\in\mathbb{N}_{0}$, then sequence
\[
\left\{ \prod_{j=1}^{m}\frac{\left(-a_{j};q_{j}\right)_{k}}{(q_{j};q_{j})_{k}}\frac{q^{\alpha k^{2}}}{{\displaystyle \prod_{r=1}^{\ell}}\left(q_{r},q_{r}^{\nu_{r}+1};q_{r}\right)_{k}}\right\} _{k=0}^{n}
\]
is PF, thus the polynomial
\[
\sum_{k=0}^{n}\prod_{j=1}^{m}\frac{\left(-a_{j};q_{j}\right)_{k}}{(q_{j};q_{j})_{k}}\frac{q^{\alpha k^{2}}z^{k}}{{\displaystyle \prod_{r=1}^{\ell}}\left(q_{r},q_{r}^{\nu_{r}+1};q_{r}\right)_{k}}
\]
has all negative zeros. For each $n\in\mathbb{N}_{0}$ and $z\in\mathbb{C}$
we have
\[
\begin{aligned} & \left|\sum_{k=0}^{n}\prod_{j=1}^{m}\frac{\left(-a_{j};q_{j}\right)_{k}}{(q_{j};q_{j})_{k}}\frac{q^{\alpha k^{2}}z^{k}}{{\displaystyle \prod_{r=1}^{\ell}}\left(q_{r},q_{r}^{\nu_{r}+1};q_{r}\right)_{k}}\right|\\
 & \le\sum_{k=0}^{\infty}\prod_{j=1}^{m}\frac{\left(-a_{j};q_{j}\right)_{k}}{(q_{j};q_{j})_{k}}\frac{q^{\alpha k^{2}}\left|z\right|^{k}}{{\displaystyle \prod_{r=1}^{\ell}}\left(q_{r},q_{r}^{\nu_{r}+1};q_{r}\right)_{k}}<\infty.
\end{aligned}
\]
Now we first apply Vitali's theorem, then apply Hurwitz's theorem
to prove that the entire function (\ref{eq:3.6}) has no zeros outside
$(-\infty,0)$. Since by Lemma 14.1.4 of \cite{Ismail1} we know that
the entire function (\ref{eq:3.6}) has infinitely many zeros, then
the assertion for on (\ref{eq:3.6}) is proved.

Let 
\[
m\ge0,\ \ell\ge1,\ \alpha=\left(\ell+\frac{m}{2}\right),\ \nu_{r}\ge0,\ 1\le r\le\ell,
\]
and 
\[
a_{j}=0,\ q_{j}=q,\ 1\le j\le m,\quad q_{r}=q,\ 1\le r\le\ell
\]
in (\ref{eq:3.6}), then for each $z\in\mathbb{N}$ and each $k\in\mathbb{N}_{0}$,
by 
\[
\lim_{q\uparrow1}\frac{(1-q)^{k}}{(q;q)_{k}}=\frac{1}{k!},\ \lim_{q\uparrow1}\frac{(1-q)^{k}}{\left(q^{\nu_{r}+1};q\right)_{k}}=\frac{1}{\left(\nu_{r}+1\right)_{k}}
\]
we have
\[
\lim_{q\uparrow1}\frac{q^{\alpha k^{2}}z^{k}\left(1-q\right)^{(m+2\ell)k}}{(q;q)_{k}^{m+\ell}{\displaystyle \prod_{r=1}^{\ell}}\left(q^{\nu_{r}+1};q\right)_{k}}=\frac{z^{k}}{\left(k!\right)^{m+\ell}{\displaystyle \prod_{r=1}^{\ell}\left(\nu_{r}+1\right)_{k}}}.
\]
Since for $n\in\mathbb{N}$ we have 
\[
nq^{n-1}\le\frac{1-q^{n}}{1-q}=1+q+\dots+q^{n-1}\le n,
\]
 and
\[
n!q^{\binom{n}{2}}\frac{(q;q)_{n}}{\left(1-q\right)^{n}}\le n!.
\]
Then, 
\[
\frac{1}{n!}\le\frac{(1-q)^{n}}{(q;q)_{n}}\le\frac{q^{-\binom{n}{2}}}{n!}<\frac{q^{-\frac{n^{2}}{2}}}{n!}.
\]
We also observe that for $\beta\ge1$ we have
\[
\frac{(1-q)^{n}}{\left(q^{\beta};q\right)_{n}}\le\frac{(1-q)^{n}}{(q;q)_{n}}.
\]
Then,
\[
\left|\frac{q^{\alpha k^{2}}z^{k}\left(1-q\right)^{(m+2\ell)k}}{(q;q)_{k}^{m+\ell}{\displaystyle \prod_{r=1}^{\ell}}\left(q^{\nu_{r}+1};q\right)_{k}}\right|\le\frac{\left|z\right|^{k}}{\left(k!\right)^{m+2\ell}}.
\]
Hence 
\[
\begin{aligned}\lim_{q\uparrow1}\sum_{k=0}^{\infty}\frac{q^{\alpha k^{2}}z^{k}\left(1-q\right)^{(m+2\ell)k}}{(q;q)_{k}^{m+\ell}{\displaystyle \prod_{r=1}^{\ell}}\left(q^{\nu_{r}+1};q\right)_{k}} & =\sum_{k=0}^{\infty}\frac{z^{k}}{\left(k!\right)^{m+\ell}{\displaystyle \prod_{r=1}^{\ell}\left(\nu_{r}+1\right)_{k}}}\end{aligned}
\]
converges uniformly for $z$ in any compact subset of $\mathbb{C}$.
Then by Hurwitz's theorem we know that the entire function 
\[
g(z)=\sum_{k=0}^{\infty}\frac{z^{k}}{\left(k!\right)^{m+\ell}{\displaystyle \prod_{r=1}^{\ell}\left(\nu_{r}+1\right)_{k}}}
\]
has no zeros outside $(-\infty,0)$. 

Let 
\[
a_{k}=\frac{1}{\left(k!\right)^{m+\ell}{\displaystyle \prod_{r=1}^{\ell}\left(\nu_{r}+1\right)_{k}}},
\]
then by the Stirling's formula \cite{Andrews1}, 
\[
\log\Gamma(x+1)=x\log x+\mathcal{O}(x),\quad x\to\infty,
\]
we have
\[
-\log a_{k}=(m+2\ell)k\log k+\mathcal{O}(k)
\]
as $k\to\infty$. Thus the order of $g(z)$ is given by \cite{Boas}
\[
\rho=\limsup_{k\to\infty}\frac{k\log k}{-\log|a_{k}|}=\frac{1}{m+2\ell},
\]
which is a positive number less than $\frac{1}{2}$. Hence $g(z)$
has infinitely many zeros, and the claim on (\ref{eq:3.7}) is proved.

\subsection{Proof of Theorem \ref{thm:func2}}

By Lemma 14.1.4 of \cite{Ismail1} we know that $_{r}A_{s}^{(\alpha)}(a_{1},\dots a_{r};\ b_{1},\dots,\ b_{s}\ ;\,q;\,z)$
has infinitely many zeros.

Let 
\[
c_{n}=\frac{(a_{1},\dots,\,a_{r};\,q)_{n}}{\left(b_{1},\dots,b_{s};\,q\right)_{n}}q^{\alpha n^{2}},
\]
then for $n\ge1$ we have 
\begin{eqnarray*}
\frac{c_{n}^{2}}{c_{n+1}c_{n-1}} & = & q^{-2\alpha}\prod_{j=1}^{r}\frac{1-a_{j}q^{n-1}}{1-a_{j}q^{n}}\cdot\prod_{k=1}^{s}\frac{1-b_{k}q^{n}}{1-b_{k}q^{n-1}}\\
 & \ge & q^{-2\alpha}\prod_{j=1}^{r}\left(1-a_{j}\right)\prod_{k=1}^{s}\left(1-b_{k}q\right)\ge4.
\end{eqnarray*}
Our assertion follows from Theorem B of \cite{Dimitrov1} or Theorem
A of \cite{Hutchinson}.


\begin{thebibliography}{10}
\bibitem{Andrews1} G. E. Andrews, R. Askey and R. Roy, \emph{Special
Functions, }Cambridge University Press, Cambridge, 1999.

\bibitem{Ahlfors1}L. Ahlfors, Complex Analysis, 3rd edition, McGraw-Hill,
1979.

\bibitem{Aissen1}M. Aissen, A. Edrei, I. J. Schoenberg, and A. M.
Whitney, On the generating functions of totally positive sequences,
Proc. Nat. Acad. Sci. U.S.A. 37 (1951), 303-307.

\bibitem{Boas}R. P. Boas, Entire Functions, Academic Press, 1954. 

\bibitem{Carnicer1}J. Carnicer, J. M. Pe\~ {n}a and A. Pinkus, On
some zero-increasing operators, Acta Math. Hungar. 94 (2002) 173-190.

\bibitem{Dimitrov1}D. K. Dimitrov and J. M. Pe\~ {n}a, Almost strict
total positivity and a class of Hurwitz polynomials, Journal of Approximation
Theory 132 (2005) 212-223.

\bibitem{Driver1} K. Driver, K. Jordaan and A. Martin\'{i}nez-Finkelshtein,
P\'{o}lya frequency sequences and real zeros of some $_{3}F_{2}$
polynomials, J. Math. Anal. 332(2007) 1045-1055.

\bibitem{Hayman}W. K. Hayman On the zeros of a q-Bessel function
Complex Analysis and Dynamical Systems, Contemp. Math., Amer. Math.
Soc., Providence, RI (2006), 205\textendash 216.

\bibitem{Hutchinson}J. I. Hutchinson, On a remarkable class of entire
functions, Trans. Amer. Math. Soc. 25 (1923) 325\textendash 332. 

\bibitem{Karlin}S. Karlin, Total Positioity, Vol. 1, Stanford University
Press, 1968.

\bibitem{Ismail1}M. E. H. Ismail, \emph{Classical and Quantum Orthogonal
Polynomials in One Variable, }Cambridge University Press, Cambridge,
2005.

\bibitem{Ismail2}M. E. H. Ismail and R. Zhang, \$q\$-Bessel Functions
and Rogers-Ramanujan Type Identities, Proceedings of American Mathematical
Society, in production.

\bibitem{Ismail3}E. H. Ismail and R. Zhang, \$q\$-Bessel Functions
and Rogers-Ramanujan Type Identities, Proceedings of American Mathematical
Society, in production.

\bibitem{Pitman}J. Pitman, Probabilistic bounds on the coefficients
of polynomials with only real zeros, Journal of Combinatorial Theory
Series. A 77 (1997) 279-303.

\bibitem{Titchmarsh1}E. C. Titchmarsh, The Theory of Functions, corrected
second edition, Oxford University Press, Oxford, 1964. \end{thebibliography}
\end{document}